\documentclass[12pt]{article}
\usepackage{latexsym,amssymb}

\def\l{{\lambda}}
\def\m{{\mu}}
\def\n{{\nu}}

\def\cR{{\cal R}}

\def\cO{{\cal O}}

\def\ZZ{\mathbb Z}
\def\CC{\mathbb C}

\def\PP{\mathbb P}

\def\proof{\noindent{\em Proof.}\enspace}

\def\qed{\nobreak\quad\nobreak$\Box$\medbreak}
\def\ot{{\mathord{\otimes }}}
\def\otc{{\mathord{\otimes\cdots\otimes }}}

\def\lra{{\mathord{\longrightarrow\;}}}

\def\Wedge{\textstyle{\bigwedge}}
\def\sg{\smallskip\noindent}
\def\mg{\medskip\noindent}
\def\bg{\bigskip\noindent}

\begin{document}

\title{Vanishing theorems for ample vector bundles}
\author{Laurent Manivel}
\date{February 1996}
\maketitle

\section{Introduction}

Since the seminal paper published by P.A. Griffiths in 1969 \cite{gr},
a whole series of vanishing theorems have been established
for the Dolbeault cohomology of ample vector bundles
on smooth projective varieties, mainly due to the efforts
of J. Le Potier, M. Schneider, A. Sommese, J-P. Demailly,
L. Ein and R. Lazarsfeld, the author, and more recently
W. Nahm \cite{de1,el,lp1,m1,m2,na,sch,so}.

This abundance of results and variants
 lead to a little confusion,
in contrast with the case of ample line bundles for which the
celebrated Kodaira-Akizuki-Nakano vanishing theorem \cite{an}
is the only general statement.

In this paper, we prove a vanishing theorem for the
Dolbeault cohomology of a product of symmetric and
skew-symmetric powers of an ample vector bundle, twisted
by a suitable power of its determinant line bundle.
This line bundle being ample, our theorem can be seen
as an effective version, in that  particular case, of the
asymptotic vanishing theorem of Serre. It includes as special
cases most of the above mentionned results.

\bg {\bf Theorem A.} {\em Let $E$ be a holomorphic
vector bundle of rank $e$,
and $L$ a line bundle on a smooth projective complex variety $X$
of dimension $n$. Suppose that $E$ is ample and $L$ nef,
or that $E$ is nef and $L$ ample.
Then, for any sequences of integers $k_1,\ldots ,k_l$ and $j_1,\ldots ,j_m$,}
$$H^{p,q}(X,S^{k_1}E\otc S^{k_l}E\ot
\Wedge^{j_1}E\otc\Wedge^{j_m}E\ot (\det E)^{l+n-p}\ot L)=0$$
{\em as soon as $p+q>n+\sum_{s=1}^m(e-j_s)$.}

\bigskip  The vanishing theorem of Griffiths, more precisely the
 version of this theorem given by M. Schneider, that is Theorem 2.4
in \cite{sch}, corresponds to $p=n$, $m=0$, $l=1$.
The extension of it  due to J-P. Demailly,
Theorem $0.2$ in \cite{de1}, corresponds to $p=n$ and $m=0$,
as will be made clear by Theorem A'.

The vanishing theorem of Le Potier, Proposition $2.3$ in  \cite{lp1},
corresponds to $p=n$, $l=0$ and $m=1$. It was generalized by L. Ein
and R. Lazarsfeld : Proposition $1.7$ of \cite{el},
at least when the vector bundles involved are identical,
corresponds to $p=n$, $l=0$ and $m$ arbitrary.

\smallskip For several ample bundles, say $E_1,\ldots ,E_h$, Theorem A
can easily be extended in the following way : the Dolbeault cohomology
groups of a product of associated bundles, each of the form appearing in
Theorem A, holds under the condition that $p+q$ be greater
than $n$ plus the sum of the contributions of each term of the
product. This does not follow from the trick used in \cite{el} to
deduce Propposition $1.7$ from Le Potier's theorem, but it is
a straightforward exercise to modify our proof of Theorem A
so as to handle the case of several vector bundles.

\smallskip Theorem A was proven in \cite{m2} for $p$ close to $n$. The methods
used  in that paper were adapted from \cite{de1} and much more complicated
than those we will use here. The main idea was to work  on flag
manifolds and compute suitable Dolbeault cohomology groups of
well-chosen line bundles -- unfortunately, there is no general
method to do that and the argument was quite involved. Here, we
restrict to relative products of projective bundles. On a projective
space, the Dolbeault cohomology of a multiple of the hyperplane
bundle is given by a special case of Bott's theorem. Actually, this
fortunate fact is a consequence of the irreducibility, as a
homogeneous bundle, of the vector bundle of regular differential forms
of a given degree on a projective space. On a grassmannian, these bundles
of forms are no longer irreducible, but remain completely reducible
--  on a more general flag manifold this property definitely disappears.

\smallskip Let us  mention that our proof of Theorem A works
{\em mutatis mutandis} in the case where
$E$ or $L$, instead of being  supposed ample, is only $k$-ample in the
sense of Sommese : one then simply has to add $k$ to the numerical
vanishing condition of the theorem to get a correct statement.
Theorem A therefore also extends Proposition $1.14$ of \cite{so}.

\medskip Let us now state our main result in a slightly different form.
Recall that if $V$ is some finite dimensional complex vector space,
the dimension of which we denote by $v$, the irreducible polynomial
representations of  the reductive algebraic group $Gl(V)$ are in
correspondance with partitions of length at most $v$, that is,
non-increasing sequences  of non-negative integers
$\l=(\l_1\geq\cdots\geq\l_v\geq 0)$. The $\l_i$ are the {\em parts} of
$\l$, while its {\em length} $l(\l)$ is the number of positive parts.
We denote the irreducible $Gl(V)$-module
corresponding to the partition $\l$ by $S_{\l}V$, and call it
the {\em Schur power} of $V$ of exponent $\l$. As special cases,
one of course recovers symmetric and skew-symmetric powers :
$$S^pV=S_{p,0,\ldots ,0}V,\quad
\Wedge^qV=S_{1,\ldots ,1,0,\ldots ,0}V,$$
where this last expression contains exactly $q$ ones.
In particular, $S_{1,\ldots ,1}V=\det V$, and one has for each
integer $m$ the factorization formula
$$S_{\l_1+m,\ldots ,\l_v+m}V=S_{\l_1,\ldots ,\l_v}V\ot (\det V)^m.$$

One usually represents a partition $\l$ by its {\em Ferrers diagram}
$D(\l)$, which is a collection of left-justified rows of decreasing lengths
from top to bottom, these lengths being given by the parts of $\l$.
If we flip this diagram over its main diagonal, we get the Ferrers
diagram of the {\em conjugate partition} $\l^*$. The parts of this
partition are simply the lengths of the columns of the Ferrers
diagram of $\l$.
We then associate to each integer $l$, with $0\leq l\leq v-1$, the integer
$$q_l(v,\l)=\sum_{j,\;\l^*_j>l}(v-\l^*_j).$$

\noindent
Since a complex vector bundle $E$ of rank $e$ on $X$ has $Gl(e,\CC)$ for
structure group, to each partition $\l$ of length at most $e$ coresponds
an associated vector bundle $S_{\l}E$.
Theorem A can then be restated in the following way :

\bg {\bf Theorem A'.} {\em Let $E$ be a vector bundle of rank $e$,
and $L$ a line bundle, on a smooth projective complex variety $X$
of dimension $n$. Suppose that $E$ is ample and $L$ nef,
or that $E$ is nef and $L$ ample. Then, for each partition $\l$
of length at most $e$, and each integer $l$ between $0$ and $e-1$, }
$$H^{p,q}(X,S_{\l}E\ot (\det E)^{l+n-p}\ot L)=0
\quad for\; p+q>n+q_l(e,\l).$$

\medskip Of course, $q_l(e,\l)$ is an decreasing function of $l$.
For each Schur power of $E$, Theorem A' therefore determines a
sufficient power of the determinant line bundle of $E$, to ensure
the vanishing of the corresponding Dolbeault cohomology in a given
degree.

\smallskip The equivalence between Theorem A and Theorem A' is an easy
consequence of Pieri's rules for the tensor product of a given Schur
power by some symmetric or skew-symmetric power. On the one hand, if
$\l$ is a partition and $l$ an integer, we let $j_a=\l_a^*$, $1\leq
a\leq m$, if $\l_a^*>l$, and $k_b=\l_b-m$ if $b\leq l$.
Then  $S_{\l}E$ is a component of
$S^{k_1}E\otc S^{k_l}E\ot \Wedge^{j_1}E\otc\Wedge^{j_m}E$,
and Theorem A thus implies Theorem A'. On the other hand,
suppose that $S_{\m}E$ is a component of $\Wedge^{j_1}E\otc\Wedge^{j_m}E$.
Its first part is then at most equal to $m$. Therefore, if $S_{\l}E$
is a component of  $S_{\m}E\ot S^{k_1}E\otc S^{k_l}E$, one has
$\l_j\leq l$ if $j>m$, so that
$$q_l(e,\l)\leq \sum_{j=1}^m(e-\l_j^*)\leq \sum_{j=1}^m(e-\m_j^*)
=\sum_{k=1}^m(e-j_k).$$
Theorem A' thus implies Theorem A.

\mg Two consequences of Theorem A-A' should be of particular interest.
First note that $q_l(e,\l)=0$ for $l\geq l(\l)$. This implies
the following partial answer to the problem raised in \cite{de1}
by J-P. Demailly :

\mg {\bf Corollary B.} {\em Under the same hypotheses,
for each partition $\l$ of length at most $e$,}
$$H^{p,q}(X,S_{\l}E\ot (\det E)^{l(\l)+n-p}\ot L)=0\quad for\;p+q>n.$$
{\em In particular, for symmetric powers,}
$$H^{p,q}(X,S^kE\ot (\det E)^{n-p+1}\ot L)=0\quad for\;p+q>n,$$
{\em while for tensor powers,}
$$H^{p,q}(X,E^{\ot k}\ot (\det E)^{\min (k,e-1)+n-p}\ot L)=0
\quad for\;p+q>n.$$

\sg Because of the factorization formula, one may always
suppose that $l(\l)<e$ in the fist part of this statement. Corollary B
therefore improves and makes more precise Theorem $0.3$ of \cite{de1}.
Moreover, its second part extends Griffiths' vanishing theorem to the
whole range of the Dolbeault cohomology. Its last part follows from the
existence of a decomposition
$$E^{\ot k}=\bigoplus_{|\l|=k}m_{\l}\; S_{\l}E,$$
where the size $|\l|$ of the partition $\l$ is defined as the sum of
its parts (the multiplicity $m_{\l}$ can be shown to be the degree of
the representation of the symmetric group, usually called a Specht
module, defined by $\l$).

\sg It was suggested in \cite{de2} that the exponent $l(\l)+n-p$
of $\det E$ in Corollary B could possibly be replaced by
$l(\l)+\min (n-p,n-q)$, and that this exponent could be the
optimal one. This was suggested by the observation that if the
Borel-Le Potier spectral sequence degenerates at the $E_2$ level,
then the exponent $e-1+\min (n-p,n-q)$ of $\det E$ is sufficient to
get a vanishing theorem in degree $p+q>n$. Unfortunately, it was
shown in \cite{m1} that this degeneracy does not occur in general.
Worse, the Borel-Le Potier spectral sequence can be non-degenerate
at an arbitrary  high level. We see no way to recover the lost
symmetry between $p$ and $q$.

\bg {\bf Corollary C.} {\em Under the same hypotheses,
for each partition $\l$ of length at most $e$,}
$$H^{p,q}(X,S_{\l}E\ot (\det E)^{n-p}\ot L)=0\quad
for\;p+q>n+q_0(e,\l),$$
{\em where $q_0(e,\l)=\sum_i(\l_1-\l_i)=e\l_1-|\l|$. In particular,}
$$H^{p,q}(X,\Wedge^kE\ot (\det E)^{n-p}\ot L)=0\quad for\;p+q>n+e-k.$$

\smallskip Sommese conjectured in \cite{so} that,
under the preceding positivity hypotheses,
$$H^{p,q}(X,\Wedge^kE\ot L)=0\quad {\rm for}\;p+q>n+e-k,$$
the case where $p=n$ being Le Potier's theorem.
This conjecture was independantly disproved in 1988 by
J-P. Demailly \cite{de1} and M. Schneider,
who gave counterexamples already for $p=n-1$.
Nevertheless,  Corollary C shows that Sommese's conjecture becomes
correct after a twist by $\det E$ to the power $n-p$.

\smallskip Apart from a few very special cases,
we have no idea about the optimality of our vanishing theorems.
Actually, it is a difficult problem to obtain interesting examples
of non-vanishing cohomology groups. To our knowledge, the only
case that has been investigated in detail is that of
grassmannians \cite{so,m2}, on which Bott's theorem allow to compute
the Dolbeault cohomology of homogeneous bundles in terms of diagrams
and the so-called {\em hook lengths}, which play an important role
in the representation theory of symmetric groups.

\sg In a final section, we propose an application of Theorem A-A'
to the geometry of degeneracy loci.
Let $\phi : E^*\lra F \ot L$ be a morphism between
vector bundles over an $n$-dimensional smooth projective variety
$X$, where $E$ and $F$ have ranks $e$ and $f$, and $L$ is a line
bundle. Denote by $D_k(\phi)$ its $k$-th degeneracy locus.
Under a strong positivity assumption, we then show that if
$D_k(\phi)$ has the expected dimension $\rho=n-(e-k)(f-k)$,
the restriction morphism $$H^q(X,{\cal O}_X)\lra H^q(D_k(\phi),
{\cal O}_{D_k(\phi)})$$ is an isomorphism for $q<\rho$, and is
injective for $q=\rho$.
A similar result holds  when $F=E$ and $\phi$ is
symmetric or skew-symmetric. In particular, $D_k(\phi)$ is
connected, if $X$ is, as soon as $\rho$ is positive.

\section{Proof of Theorem A}
\subsection{On the Borel-Le Potier spectral sequence}

Let again $E$ be a complex vector bundle of rank $e$, $F$ another
vector bundle,  on a smooth compact complex
variety $X$ of dimension $n$. Let $Y=\PP(E^*)$ be the variety of
hyperplanes of $E$, let $\pi$ be its projection onto $X$. We denote as
usual by  $\cO_E(1)$ the universal quotient line bundle on $Y$, and by
$\cO_E(k)$ its $k$-th power, $k\in\ZZ$.

To compute the Dolbeault cohomology of $\cO_E(k)\ot\pi^*F$
on $Y$, we use the fact that the bundle $\Omega_Y^p$
of regular $p$-forms on $Y$ may be filtered according to their
degree $t$ on $X$, that is, by the subbundles
$$F^{t,p}=\pi^*\Omega_X^t\Wedge \Omega_Y^{p-t}.$$
The corresponding quotient are the vector bundles
$$G^{t,p}=F^{t,p}/F^{t+1,p}=\pi^*\Omega_X^t\ot \Omega_{Y/X}^{p-t},$$
where $\Omega_{Y/X}^m$ denotes the vector bundle of relative
differential forms of degree $m$ -- so that its restriction to
each fiber of $\pi$ is the bundle of $m$-forms on that fiber.

 This implies, for each integer $p$,
the existence of a spectral sequence, which we call after \cite{de1} the
{\em Borel-Le Potier spectral sequence} \cite{lps1}.
Its  $E_1$-term is given by
$${}^pE_1^{t,s}=H^{t+s}(Y,\Omega_{Y/X}^{p-t}
\ot \cO_E(k)\ot\pi^*(\Omega_X^t\ot F)),$$
 The Borel-Le Potier spectral sequence converges
to the Dolbeault cohomology groups $H^{p,q}(Y,\cO_E(k)\ot\pi^*F)$,
that is, there is a naturally defined filtration of that complex
vector space, with associated graded quotient
$${\rm gr}\;H^{p,q}(Y,\cO_E(k)\ot\pi^*F)=
\bigoplus_{t+s=q}{}^pE_{\infty}^{t,s}.$$

Let us compute the $E_1$-terms of the Borel-Le Potier spectral
sequence. We will denote by $(k|l)$, $k\geq 0$, $l>0$, the partition
of length $l$,  whose first part is equal to $k+1$, while its other
non zero parts are equal to one : its Ferrers diagram is a hook.
The corresponding bundle $S_{(k|l)}E$ associated to $E$ is
the kernel of the natural contraction map
$$S^{k+1}E\ot \Wedge^{l-1}E\lra S^{k+2}E\ot E^*\ot\Wedge^{l-1}E
\lra S^{k+2}E\ot \Wedge^{l-2}E.$$
     In particular,     $S_{(k|1)}E=S^{k-1}E$, $S_{(k|e)}E=S^{k}E
\ot \det E$,     $S_{(0|j)}E=\Wedge^jE$.
We keep the same notation when $k<0$ or $l>e$,
in which case $S_{(k|l)}E=0$   (except for $k=-1$ and $l=0$).

\mg {\bf Lemma D.} {\em  Suppose that $k>0$. Then }
$${ }^pE_1^{t,s-t}=H^{t,s}(X,S_{(k-p+t-1|p-t+1)}E\otimes F).$$

\proof We use a Leray spectral sequence to compute the cohomology group
$${}^pE_1^{t,s-t}=H^{  s}(Y,\Omega_{Y/X}^{p-t}
\ot \cO_E(k)\ot \pi^*(\Omega_X^t\ot F)).$$
Let $T$ denote the tautological hyperplane bundle, of rank $e-1$ on $Y$.
The bundle of relative one-forms is $\Omega_{Y/X}^1=\cO_E(-1)\ot T$,
so that its wedge powers $\Omega_{Y/X}^h=\cO_E(-h)\ot \Wedge^hT$.
Such a bundle is homogenous and irreducible, and its cohomology
is therefore  given by Bott's theorem (\cite{dem}, see also \cite{lp2}) :
$$\cR^i_{\pi *}(\Omega_{Y/X}^{p-t}\ot \cO_E(k))=
\cR^i_{\pi *}(\cO_E(k-p+t)\ot\Wedge^{p-t}T)=
\delta_{i,0}\; S_{(k-p+t-1|p-t+1)}E.$$
The Leray spectral sequence then degenerates at $E_2$,
and the lemma follows. \qed

\mg More generally, let $Y_h$ be the relative product over $X$, of $h$
copies of $Y$, and denote the corresponding product of
powers  of exponents $k_1,\ldots ,k_h$ of the universal quotient
line bundle by $\cO_E(k_1,\ldots ,k_h)$. The Dolbeault cohomology of
this line bundle, twisted by the pull-back of $F$, is again
computed by a Borel-Le Potier spectral sequence, whose $E_1$ term
is given by the following straightforward generalization of
lemma D \nolinebreak :

\mg {\bf Lemma D'.} {\em  Suppose that $k_1,\ldots ,k_h>0$. Then}
$${ }^pE_1^{t,s-t}=\bigoplus_{p_1+\ldots +p_h=p-t}
H^{t,s}(X,S_{(k_1-p_1-1|p_1+1)}E\otc S_{(k_h-p_h-1|p_h+1)}E \ot F).$$

\sg Such a cohomology group can be non zero only if
$p_i<\min (e,k_i)$ for all $i\leq h$, hence
$${ }^pE_1^{t,s-t}=0\quad
{\rm if}\;t<p-\sum_{i=1}^h\min (e-1,k_i-1).$$

\subsection{The induction}

We now proceed by induction on $r=n-p$.
For simplicity, we divide it in two steps.
Here, $E$ is a vector bundle of rank $e$, and $L$ a line bundle
on $X$. Moreover, we make the positivity hypothesis that $E$ is nef
and $L$ ample, or that $E$ is ample and $L$ nef.

\mg {\bf First step.} Let us prove that for any $r\geq 0$, Theorem A
for $n-p<r$ implies Theorem A for $n-p=r$ and $l=0$ -- that is,
Corollary C for $n-p=r$.
A collection of integers $j_1,\ldots ,j_m$ being given,
with $j_k\leq e$ for all $k$, let us consider the Borel-Le Potier
spectral sequence that computes the Dolbeault cohomology of
$\cO_E(j_1,\ldots ,j_m)\ot\pi^*((\det E)^r\ot L)$ on $Y_m$.
If $p_0=n-r+\sum_{k=1}^m(j_k-1)$, our previous lemma implies that
$${ }^{p_0}E_1^{n-r,s-n+r}=H^{n-r,s}(X,
\Wedge^{j_1}E\otc\Wedge^{j_m}E\ot (\det E)^r\ot L),$$
which is the cohomology group we want to show to be zero.
Moreover, for each positive integer $u$,
${ }^{p_0}E_1^{n-r-u,s-n+r+u}=0$, while ${ }^{p_0}E_1^{n-r+u,s-n+r-u}$
is given by the direct sum
$$\bigoplus_{u_1+\ldots +u_m=u} H^{n-r+u,s}
(X,S_{(u_1|j_1-u_1)}E\otc S_{(u_m|j_m-u_m)}E\ot (\det E)^r\ot L).$$
To show that such a            group vanishes, we observe that
for each $k$, $S_{(u_k|j_k-u_k)}E$
is a direct factor of $S^{u_k}E\ot\Wedge^{j_k-u_k}E$, so that
it suffices to prove that
$$H^{n-r+u,s} (X,S^{u_1}E\otc S^{u_m}E\ot
\Wedge^{j_1-u_1}E\otc \Wedge^{j_m-u_m}E\ot (\det E)^r\ot L) $$
vanishes. But at most
$u$ of the above    symmetric powers can have positive exponents,
so that we have enough copies of $\det E$ to apply our
induction hypothesis in degree $n-r+u$.  We obtain
$${ }^{p_0}E_1^{n-r+u,s-n+r-u}=0
\quad {\rm for}\;s>r+\sum_{k=1}^m(e-j_k).$$
    We then consider, for $u>0$, the following morphisms
 of our                 spectral sequence :
$${ }^{p_0}E_u^{n-r-u,s-n+r-1+u}\lra { }^{p_0}E_u^{n-r,s-n+r}
\lra { }^{p_0}E_u^{n-r+u,s-n+r+1-u}.$$
The left hand side is zero because ${ }^{p_0}E_1^{n-r-u,s-n+r-1+u}$
is  zero, and the right hand side vanishes, as      we have just
seen,  as soon as $s\geq r+\sum_{k=1}^m(e-j_k)$. Hence       an
isomorphism
$${ }^{p_0}E_1^{n-r,s-n+r}\simeq { }^{p_0}E_{\infty}^{n-r,s-n+r}.$$
But this last term is a graded piece of some filtration of
$$H^{p_0,s}(Y_m,\cO_E(j_1,\ldots,j_m)\ot\pi^*((\det E)^r\ot L)).$$
Since our positivity assumption on $E$ and $L$, and the fact
that $j_1,\ldots, j_m$ are all positive,
implies that the line bundle
$\cO_E(j_1,\ldots ,j_m)\ot \pi^*((\det E)^r\ot L)$ is ample,
this Dolbeault cohomology group vanishes, by the
Kodaira-Akizuki-Nakano vanishing theorem, as soon as $p_0+s>\dim Y_m$,
that is, $s>r+\sum_{k=1}^m(e-j_k)$. This is precisely what we wanted to prove.

\mg {\bf Second step.} We now show that Theorem A follows for $n-p=r$
and $l$ arbitrary.
We fix two collections of integers $k_1,\ldots   ,k_l$ and
$j_1,\ldots   ,j_m$, and consider the Borel-Le Potier spectral sequence
that computes the Dolbeault cohomology of $$\cO_E(k_1+e,\ldots ,k_l+e)
\ot\pi^*(\Wedge^{j_1}E\otc \Wedge^{j_m}E\ot (\det E)^r\ot L)$$
on $Y_l$.
If $p_1=n-r+l(e-1)=\dim Y_l-r$, we have, again by the previous lemma,
$$^{p_1}E_1^{n-r,s-n+r}=H^{n-r,s}(X,S^{k_1}E\otc S^{k_l}E\ot
\Wedge^{j_1}E\otc\Wedge^{j_m}E\ot (\det E)^{l+r}\ot L),$$
and this is the group we want to show to be zero. Moreover, for each
positive integer $v$, we have ${ }^{p_1}E_1^{n-r-v,s-n+r+v}=0$, while
${ }^{p_1}E_1^{n-r+v,s-n+r-v}$ is the direct sum
$$\begin{array}{l}
\bigoplus_{v_1+\ldots +v_l=v}
H^{n-r+v,s}(X,S_{(k_1+v_1|e-v_1)}E\otc S_{(k_l+v_l|e-v_l)}E\ot \\
\hspace*{5cm}\ot\Wedge^{j_1}E\otc\Wedge^{j_m}E\ot (\det E)^{l+r}\ot
L).
\end{array} $$
But now, for each  $i$, $S_{(k_i+v_i|e-v_i)}E$
is a direct factor of $S^{k_i+v_i}E\ot \Wedge^{e-v_i}E$, and at least
$l-v$ of these wedge powers are copies of $\det E$. We therefore have
enough copies of the determinant line bundle to use our induction
hypothesis  in degree $n-r+v$, and we obtain
$${ }^{p_1}E_1^{n-r+v,s-n+r-v}=0
\quad {\rm   for}\;s>r+\sum_{k=1}^m(e-j_k).$$
Then we consider, for  $v>0$, the following morphisms
 of our                 spectral sequence :
$${ }^{p_1}E_v^{n-r-v,s-n+r-1+v}\lra { }^{p_1}E_v^{n-r,s-n+r}
\lra { }^{p_1}E_v^{n-r+v,s-n+r+1-v}.$$
The left hand side is zero, because ${ }^{p_1}E_1^{n-r-v,s-n+r-1+v}$
is        zero, and the right hand side vanishes, as      we have just
seen,  for        $s\geq r+\sum_{k=1}^m(e-j_k)$. Hence            an
isomorphism
$${ }^{p_1}E_1^{n-r,s-n+r}\simeq { }^{p_1}E_{\infty}^{n-r,s-n+r}.$$
But this last term is a graded piece of some filtration of
   $$H^{p_1,s}(Y_l,\Wedge^{j_1}\pi^*E\otc \Wedge^{j_m}\pi^*E\ot (\det
\pi^*E)^r  \ot\cO_E(k_1+e ,\ldots,k_l+e)\ot\pi^*L ).$$
Since $\pi^*E$ is nef and $\cO_E(k_1+e ,\ldots,k_l+e)\ot\pi^*L$
is ample on $Y_l$,         our first step implies the vanishing
of the previous cohomology group for $s>r+\sum_{k=1}^m(e-j_k)$.
Theorem A is proved.

\section{An application to degeneracy loci}

It was suggested in \cite{lay} that vanishing theorems could
provide interesting information on the geometry of degeneracy loci.
We shall illustrate this idea in that section. It mainly consists
in making use of a resolution of the ideal sheaf of a given
degeneracy locus. Then we apply vanishing theorems to the
different vector bundles involved in this resolution, so as to
get some control of the cohomology of that ideal. This allows
for example to prove the connectedness of such a degeneracy loci,
although under a rather strong positivity hypotheses,
and to go a little further.

The simplest case is that of a morphism $\phi : E^*\lra F$
between vector bundles of ranks $e$ and $f$, on a smooth projective
variety $X$ of dimension $n$. The $k$-th degeneracy locus
$$D_k(\phi)=\{x\in X,\;{\rm rank}\;\phi_x\leq k\}$$
has a natural scheme structure, given in a trivialisation by the
vanishing of  minors of  order $k+1$. We make the
assumption that $D_k(\phi)$ has the {\em expected codimension},
namely $(e-k)(f-k)$. Its ideal sheaf ${\cal I}$ then has
Lascoux's complex $K^{\bullet}\lra {\cal I}\lra 0$
for minimal free resolution \cite{las}, where
$$K^i=\bigoplus_{|\l|=i}S_{\l(k)}E^*\ot S_{\l^*(k)}F^*.$$
The notation $\l(k)$ in that expression means the following : if $l$
is the {\em rank} of $\l$, that is the side of the largest square
contained in its diagram, then the partition $\l(k)$ is obtained
by adjoining $k$ parts equal to $l$ to $\l$. To get a nonzero Schur
power, we therefore need that $\l_1+k\leq f$ and $\l_1^*+k\leq e$.
We will write $\l=(l,\m,\n)$ if $\m$ and $\n$ are the partitions
defined by $\m_i=\l_i-l$, $1\leq i\leq l$, and $\n_j=\l_{j+l}$.

\bg {\bf Theorem E.} {\em Let $\phi : E^*\lra F\ot L$ be a morphism
between vector
bundles on a smooth projective variety $X$ of dimension $n$, where $E$
has rank $e$, $F$ has rank $f$, and $L$ is a line bundle. Let $k<
\min (e,f)$ be a positive integer. We make the following assumptions :
\begin{itemize}
\item $E$ is ample and $F$ is nef, or $E$ is nef and $F$ ample,
\item $L^{\ot k}\geq \det E\ot\det F$,
\item $D_k(\phi)$ has the expected dimension $\rho=n-(e-k)(f-k)$.
\end{itemize}
Then the natural restriction map $H^q(X,{\cal O}_X)\lra
H^q(D_k(\phi),{\cal O}_{D_k(\phi)})$ is an isomorphism for
$0\leq q<\rho$, and is still injective for $q=\rho$.}

\bigskip If $A$ and $B$ are line bundles, we mean by $A\geq B$ that
$A\ot B^*$ is nef.

\mg {\em Remark.}
If $X$ is connected, $D_k(\phi)$ is connected as soon as
$\rho>0$ : this was proved in \cite{fl} for $L={\cal O}_X$ and $E\to
F$ ample, without any of our extra hypotheses. But to our
knowledge, the equality
$h^{0,q}(D_k(\phi))=h^{0,q}(X)$ for $q<\rho$
is new. Note that there is in general no Barth-Lefschetz type
isomorphism theorem for degeneracy loci (think for example to the
Segre embedding of $\PP^m\times \PP^n$), although the map
$$H_q(D_k(\phi),\ZZ)\lra H_q(X,\ZZ)$$
is always surjective for $q\leq \rho$ \cite{fl}.

\mg\proof Let again ${\cal I}$ be the ideal sheaf of $D_k(\phi)$.
Our statement is equivalent to the vanishing of $H^q(X,{\cal I})$
for $q\leq\rho$, so that it is enough to verify that
$$H^{q+i-1}(X,K^i)=0\quad {\rm for}\;0<q\leq\rho,\;i>0.$$
Now, by Serre duality, this is equivalent to the vanishing
$$H^{n,n-q-|\l|+1}(X,S_{\l(k)}E\ot S_{\l^*(k)}F\ot L^{\ot |\l(k)|})=0$$
for $q\leq\rho$ and each non-empty partition $\l$. But if $\l$ has
rank $l$, then $|\l(k)|\geq kl$ and $L^{\ot |\l(k)|}\geq
(\det E)^{\ot l}\ot (\det F)^{\ot l}$. Theorem A', and the remarks
following Theorem A, therefore imply that
$$H^{n,p}(X,S_{\l(k)}E\ot S_{\l^*(k)}F\ot L^{\ot |\l(k)|})=0$$
for $p>q_l(\l(k))=l(e-k-l)-|\n|+l(f-k-l)-|\m|$, with the notations of
the previous picture. But for $p=n-q-|\l|+1$, this condition is
equivalent to $$q\leq n-(e-k)(f-k)+(e-k-l)(f-k-l),$$
and the rightmost term of that inequality is non negative. \qed

\mg The preceeding approach is also suited to the case of a
morphism $\phi : E^*\lra E\ot L$ which is supposed to be symmetric
or skew-symmetric. The expected dimension of the degeneracy locus
$D_k(\phi)$ is $\rho=n-{e-k+1\choose 2}$ in the symmetric case, and
$\rho=n-{e-k\choose 2}$ in the skew-symmetric case.
When this dimension is correct, the minimal resolution $K^{\bullet}
\lra {\cal I}\lra 0$ of the ideal sheaf of $D_k(\phi)$ has been
computed in \cite{jpw} : if $\phi$ is symmetric,
$$K^i=\bigoplus_{\l=(l,\m,\m^*),\;i=|\m|+l(l-1)/2}S_{\l(k-1)}E^*,$$
where $l$ has to be even, while if $\phi$ is skew-symmetric and $k$ is even,
$$K^i=\bigoplus_{\l=(l,\m,\m^*),\;i=|\m|+l(l+1)/2}S_{\l(k+1)}E^*.$$
With the very same proof, Theorem E extends in the
following way :

\bg {\bf Theorem F.} {\em Let $\phi : E^*\lra E\ot L$ be a symmetric or
  skew-symmetric morphism on a smooth projective variety $X$ of
dimension $n$, where $E$
has rank $e$ and $L$ is a line bundle. Let $k<e$ be a positive
integer, with $k$ even if $\phi$ is skew-symmetric.
We make the following assumptions :
\begin{itemize}
\item $E$ is ample, or simply nef if $L$ is ample,
\item $L^{\ot k}\geq \det E$,
\item $D_k(\phi)$ has the expected dimension $\rho=n-{e-k+1\choose 2}$
  if $\phi$ is symmetric, and $\rho=n-{e-k\choose 2}$ if $\phi$ is
  skew-symmetric.
\end{itemize}
Then the natural restriction map
$H^q(X,{\cal O}_X)\lra H^q(D_k(\phi),{\cal O}_{D_k(\phi)})$
is an isomorphism for $0\leq q<\rho$, and is injective for $q=\rho$.}

\bg {\bf Corollary G.} {\em Under the same hypotheses, $D_k(\phi)$ is
connected as soon as $\rho>0$.}

\mg This connectedness was established in \cite{tu}, without our extra
hypotheses, in the skew-symmetric case, and in the symmetric case but
only when the rank  $e$ is even. Again we obtain extra information when the
dimension of $D_k(\phi)$ is greater than one.

It would be interesting to know whether the conclusions of Theorems E
and F hold true for a morphism
 $\phi : E^*\lra F$ under the only assumption that
$E\ot F$ is ample for the former, and $\phi : E^*\lra E\ot L$ with
$S^2E\ot L$ or $\Wedge^2E\ot L$ ample for the latter, at least when the
degeneracy locus has the expected dimension. We hope that suitable
vanishing theorems will lead to this improvement of Fulton, Lazarsfeld
and Tu's results.

\nopagebreak\mg Institut Fourier, UMR 5582 UJF/CNRS,
Universit\'e de Grenoble I, BP 74, 38402 Saint Martin d'H\`eres,
France. e-mail : manivel@puccini.ujf-grenoble.fr

\end{document}